\documentclass[a4paper,12pt]{article} 
\newdimen\paperhight
\usepackage{amsmath}
\usepackage{amssymb}
\usepackage{amsfonts}
\usepackage[usenames]{color}

\newcommand{\pd}{\partial} 

\newcommand{\pr}{\par \vspace{3mm}\noindent [{\bf Proof}] \qquad}
\newcommand{\prend}{\hfill \qed \par \vspace{3mm}}

\newcommand{\qed}{\quad\hbox{\rule[-2pt]{3pt}{6pt}}\par\vspace{3mm}}

\newcommand{\1}{{\bf 1}} 
 
\newcommand{\C}{\mathbb C} 
\newcommand{\Z}{\mathbb Z}

\newcommand{\N}{\mathbb N}

\newcommand{\CI}{{\cal I}}

\newcommand{\CL}{{\cal L}}

\newcommand{\CU}{{\cal U}}

\newcommand{\CY}{{\cal Y}}

\newcommand{\CH}{{\cal H}}

\newcommand{\tr}{{\rm tr}}
\newcommand{\wt}{{\rm wt}}
\newcommand{\Hom}{{\rm Hom}}

\newcommand{\Image}{{\rm Im}}
\newcommand{\Ker}{{\rm Ker}}

\newcommand{\Tr}{{\rm Tr}}
\newtheorem{thm}{Theorem}
\newtheorem{prn}[thm]{Proposition}

\newtheorem{lmm}[thm]{Lemma}

\begin{document}
\title{Flatness of Tensor Products and Semi-Rigidity for $C_2$-cofinite Vertex Operator Algebras II}
\author{\begin{tabular}{c}
Masahiko Miyamoto\\
Institute of Mathematics, \\
University of Tsukuba, \\
Tsukuba, 305 Japan \end{tabular}}
\date{}
\maketitle

\begin{abstract}
Let $V$ be a simple $C_2$-cofinite VOA of CFT type and we study 
the properties of non-semi-simple modules. We assume that 
there is a $V$-module $Q$ such that $\Hom_V(Q\boxtimes V',V)\not=0$.  
Let us consider a trace function $\Psi_{V}^{\tr}$ on $V$. 
As the author has shown in \cite{M1}, 
an $S$-transformation $S(\Psi_V^{\tr})$ 
of $\Psi_{V}^{\tr}$ corresponding to 
$\begin{pmatrix}0&-1\cr 1&0\end{pmatrix}$ may contain 
pseudo-trace functions. 
In this paper, we assume that $S(\Psi_V^{\tr})$ is a linear combination 
of trace functions (i.e. no pseudo-trace functions), then 
we show that all $V$-modules are semi-rigid and all trace functions 
$\Psi_{U}^{\tr}$ on simple modules $U$ appear in $S(\Psi_V^{\tr})$.  
As such an example, we show that a $C_2$-cofinite orbifold model $V$ of a rational 
VOA of CFT type  
has no pseudo-trace functions in $S(\Psi_V^{\tr})$. 
As a corollary of our main theorem, such an orbifold model becomes rational.   
\end{abstract} 

\section{Introduction}
Let $V=\oplus_{n=0}^{\infty}V_n$ be a simple $C_2$-cofinite 
vertex operator algebra of CFT type (i.e. $\dim V_0=1$). 
Since $V$ is $C_2$-cofinite, a fusion product $W\boxtimes U$ of 
finitely generated $V$-modules $W$ and $U$ 
is well-defined as a maximal finitely generated $V$-module with a 
surjective intertwining operator of $W$ from $U$ to $W\boxtimes U$. 
From the maximality, we can induced a canonical homomorphism 
$\delta\boxtimes {\rm id}_U:W\boxtimes U\to W^1\boxtimes U$ from 
a $V$-homomorphism $\delta:W\to W^1$, 
where ${\rm id}_U$ denotes the identity map on $U$.  We assume: \\

\noindent
Hypothesis I: \quad For every irreducible $V$-module $W$, there is 
an irreducible $V$-module $\widetilde{W}$ such that $\Hom_V(W\boxtimes \widetilde{W},V)\not=0$. \\

Since $V$ is $C_2$-cofinite, the associativity of fusion products holds and 
so Hypothesis I is equivalent to the existence of $\widetilde{V'}$ for 
a restricted dual $V'$ of $V$ because of  
$$(\widetilde{V'}\boxtimes W')\boxtimes W\cong \widetilde{V'}\boxtimes (W\boxtimes W')\xrightarrow{\mbox{epi}} \widetilde{V'}\boxtimes V' \xrightarrow{\mbox{epi}}V.$$ 
In the case $V'\cong V$, then we can take a restricted dual 
$W'$ of $W$ as $\widetilde{W}$. 

The author introduced a concept of "Semi-Rigidity" in \cite{M2}. 
We call an irreducible $V$-module $W$ {\bf semi-rigid} 
if there are epimorphisms $e_W:W\boxtimes \widetilde{W}\to V$ and 
$e_{\widetilde{W}}:\widetilde{W}\boxtimes W\to V$ and 
a homomorphism $\rho:P \to \widetilde{W}\boxtimes W$ 
satisfying $e_{\widetilde{W}}\rho(P)=V$ such that 
$$(e_W\boxtimes {\rm id}_W)(\mu({\rm id}_W\boxtimes \rho)(W\boxtimes P))=V\boxtimes W,$$ 
where 
$\mu:W\boxtimes (\widetilde{W}\boxtimes W) \rightarrow 
(W\boxtimes \widetilde{W})\boxtimes W$ 
is a canonical isomorphism (see (2.2)), and $P$ is a projective cover of $V$. 
Namely, we consider the following diagram
$$\begin{array}{ccccc}
W\boxtimes \rho(P)&\subseteq &W\boxtimes (\widetilde{W}\boxtimes W)&\xrightarrow{\mu} &(W\boxtimes \widetilde{W})\boxtimes W \cr
\mbox{}\qquad\qquad  \downarrow {\rm id}_W\boxtimes e_{\widetilde{W}}&          
&\downarrow    && \mbox{}\qquad\qquad\downarrow e_{W}\boxtimes {\rm id}_W\cr
W\boxtimes V&\cong     &W\boxtimes V                         && V\boxtimes W
\end{array}$$ 
We call $V$ semi-rigid when all simple $V$-modules are semi-rigid.  

As the author showed in \cite{M1}, if we consider a trace function 
$\Psi_{V}^{\tr}(\ast,\tau)$ on $V$ by 
$$\Psi_{V}^{\tr}(v,\tau)=\Tr_V (o(v) q^{\tau(L(0)-c/24)})
=\sum_{n=0}^{\infty}(\Tr_{V_n}o(v))q^{(n-c/24)\tau} \quad 
\mbox{ for }v\in V,$$
then its $S$-transformation $S(\Psi_{V,\tr})$ (corresponding to 
$\begin{pmatrix}0&-1\cr 1&0\end{pmatrix}\in SL_2(\Z)$), which is given by 
$(\frac{-1}{\tau})^{\wt(v)}\Psi_V^{\tr}(v,-1/\tau)$, equals to a linear combination of 
trace functions $\Psi_{U}^{\tr}$ on simple modules $U$ and 
pseudo-trace functions $\Psi_T^{\phi}$,  
where $q^{\tau}$ denotes $e^{2\pi i\tau}$, $c$ is a central charge of $V$ and 
$o(v)$ denotes a grade-preserving operator of 
$v$, (e.g. $o(v)=v_{m-1}$ for $v\in V_m$). In other words, there are 
$\lambda_{(U,\tr)}, \lambda_{(T,\phi)}\in \C$ such that for 
$v\in V$ and $0<|q^{\tau}|<1$, we have:
$$S(\Psi_V^{\tr})(v,\tau)=\sum_{U\mbox{ irr.}} \lambda_{(U,\tr)}\Psi_{U}^{\tr}(v,\tau)
+\sum_{(T,\phi),\phi\not=\tr} \lambda_{(T,\phi)}\Psi_T^{\phi}(v,\tau).$$

\noindent
Hypothesis II: \quad $S(\Psi_V^{\tr})$ is a linear combination of trace functions on modules. \\

The aim of this paper is to show the following theorem. \\

\noindent
{\bf Main Theorem} \quad {\it
If $V$ is a simple $C_2$-cofinite vertex operator algebra of CFT-type 
satisfying Hypothesis I and II, 
then $V$ is semi-rigid and $\lambda_{(U,\tr)}\not=0$ for every simple module $U$. }\\

At last, we will show an example satisfying Hypothesis II. \\

\noindent
{\bf Theorem 4} \quad {\it  
Let $T$ be a rational vertex operator algebra of CFT type and 
$\tau$ is a finite automorphism of $T$.  We assume that 
the fixed point subVOA $T^{\tau}$ is $C_2$-cofinite and satisfies Hypothesis I. 
Then $T^{\tau}$ has no pseudo-trace functions in $S(\Psi_{T^{\tau}}^{\tr})$. }\\ 

Here a VOA $T$ is called rational if all $\N$-gradable modules are 
completely reducible. As a corollary, we have:

\noindent
{\bf Corollary 7} \quad {\it  
Under the assumption of Theorem 4, $T^{\tau}$ is also rational.}\\

\noindent
{\rm Acknowledgement} \\
The author wishes to thank T.~Abe and H.~Yamauchi for their pointed questions.   
He also thanks to G.~H\"{o}hn, V.~Nikulin and N.~Scheithauer, the organizers of 
the conference held at Edinburgh in 2009 where he was stimulated to extend 
his result from the assumption $V\cong V'$ to general cases.

\section{Preliminary results}
\subsection{The space of logarithmic intertwining operators}
Throughout this paper, we will treat only finitely generated modules and 
${\rm mod}(V)$ denotes the set of all finitely generated $V$-modules. 
We should note that all modules are $\N$-gradable because $V$ is $C_2$-cofinite. 
If $V$ is also rational, then 
Proposition 5.1 in \cite{HV} implies that $V$ is rigid (which is stronger than semi-rigid). 
Our aim is to extend his result to non-semi-simple modules. 
Therefore, our concept of intertwining operators includes 
logarithmic intertwining operators $\CY$ of type $\binom{A}{U,B}$ which has a form 
$\CY(u,z)=\sum_{i=0}^K\sum_r u_{r,i}z^{-r-1}\log^iz$ with $u_{r,i}\in \Hom(B,A)$ for 
$u\in U$. 
$Y^U$ denotes a vertex operator of $V$ on a module $U$,  
and $\CI_{A,B}^C$ denotes the space of (logarithmic) intertwining operators 
of type $\binom{C}{A\quad B}$ for $A,B,C\in {\rm mod}(V)$. 
We fix one surjective intertwining operator 
$\CY_{A,B}^{\boxtimes}\in \CI_{A,B}^{A\boxtimes B}$ for each pair $A,B$. 
Here "surjective" implies 
that the images of $\CY_{A,B}^{\boxtimes}$ are not contained a proper subspace. 
For $W\in {\rm mod}(V)$, $W_r$ denotes a generalized eigenspace 
$\{w\in W \mid (L(0)-r)^Nw=0 \quad \mbox{ for some }N\in \N\}$ and 
$W'$ denotes a restricted dual $V$-module $\oplus_r \Hom(W_r,\C)$ of $W$.

\subsection{Analytic functions} 
We first recall the analytic part on the 
composition of intertwining operators (with logarithmic terms) from \cite{H1}. 
From now on, let $A,B,C,D,E,F\in {\rm mod}(V)$ 
and $a\in A, b\in B, c\in C, d'\in D'$. As Huang showed,  
for intertwining operators $\CY_1\in \CI_{A,E}^D$, $\CY_2\in \CI_{B,C}^E$, 
$\CY_3\in \CI_{F,C}^D$ and $\CY_4\in \CI_{A,B}^F$, 
the formal power series (with logarithmic terms) 
$$ \begin{array}{ccc}
\langle  d',\CY_1(a,x)\CY_2(b,y)c\rangle  & \mbox{ and } &
\langle  d',\CY_3(\CY_4(a,x-y)b,y)c\rangle 
\end{array}$$
are absolutely convergent in  $\Delta_1=\{(x,y)\in \C^2\mid|x|>|y|>0\}$ 
and $\Delta_2=\{(x,y)\in \C^2 \mid |y|>|x-y|>0\}$, respectively, 
and can all be analytically extended to multi-valued analytic functions on 
$$M^2=\{(x,y)\in \C^2 \mid xy(x-y)\not=0 \}.$$
As he did, we are able to lift them to single-valued analytic functions 
$$\begin{array}{ccc}
E(\langle  d,\CY_1(a,x)\CY_2(b,y)c\rangle ) &\mbox{ and }&
E(\langle  d,\CY_3(\CY_4(a,x-y)b,y)c\rangle ) 
\end{array}$$
on the universal covering $\widetilde{M^2}$ of $M^2$. 
As he remarked, single-valued liftings are not unique, but 
the existence of such functions is enough for our arguments. 
The important fact is that if we fix $A,B,C,D$, then these functions 
are given as solutions of the same differential equations. 
Therefore, for $\CY_1\in \CI_{A,E}^D,\CY_2\in \CI_{B,C}^E$ there are $\CY_5\in \CI_{A\boxtimes B,C}^D$ and 
$\CY_6\in \CI_{B,A\boxtimes C}^D$ 
such that 
$$\begin{array}{cl}
E(\langle d',\CY_1(a,x)\CY_2(b,y)c\rangle)
=E(\langle d',\CY_5(\CY_{A,B}^{\boxtimes}(a,x-y)b,y)c\rangle) &\mbox{ and}\cr
E(\langle d',\CY_2(\CY_4(a,x-y)b,y)c\rangle)
=E(\langle d',\CY_6(a,x)\CY_{B,C}^{\boxtimes}(b,y)c\rangle). &
\end{array}\eqno{(2.1)}$$ 
We note that the right hand sides of (2.1) are usually expressed by linear sums, say, 
$$E(\langle d',\CY_1(a,x)\CY_2(b,y)c\rangle)=\sum_i E(\langle d',\CY_{1i}(\CY_{2i}(a,x-y)b,y)c\rangle).  $$
For each term, from the maximality of fusion products, 
there is a homomorphism $\xi_i\in \Hom_V(A\boxtimes B, \Image(\CY_{2i}))$ such that 
$\CY_{2i}=\xi_i\CY_{A,B}^{\boxtimes}$. Then it is easy to check that 
$\sum_i \CY_{1i}\xi_i$ is an intertwining operator in 
$\CI_{A\boxtimes B,C}^D$ and so we can get the expressions (2.1).  
The canonical isomorphism $\mu: (A\boxtimes B)\boxtimes C \to 
A\boxtimes(B\boxtimes)$ is 
given by 
$$E(\langle d',\mu
\CY_{A\boxtimes B, C}^{\boxtimes}(\CY_{A,B}^{\boxtimes}(a,x-y)b,y)c\rangle)
=E(\langle d',\CY_{A,B\boxtimes C}^{\boxtimes}(a,x)\CY_{B,C}^{\boxtimes}(b,y)c\rangle). 
\eqno{(2.2)}$$

\subsection{Skew symmetric and adjoint intertwining operators}
In his paper \cite{HV}, Huang explicitly defined a skew symmetry intertwining operator 
$\sigma_{12}(\CY)\in \CI_{B,A}^C$ and an adjoint intertwining operator 
$\sigma_{23}(\CY)\in \CI_{A,C'}^{B'}$ for $\CY\in \CI_{A,B}^C$ 
under the assumption that $\CY$ has no logarithmic terms.  
Even if $\CY\in \CI_{B,A}^C$ has logarithmic terms, 
by considering a path $\{z=\frac{1}{2}e^{\pi i t}x \mid t\in [0,1]\}$, 
there is $\widetilde{\CY}\in \CI_{A,B}^C$ such that 
$$E(\langle c', \widetilde{\CY}(a,z)\sigma_{12}(Y^B)(b,x)\1\rangle)
=E(\langle c',\CY(b,x)\sigma_{12}(Y^A)(a,z)\1\rangle), \eqno{(2.3)}$$
which implies there is an isomorphism $\sigma_{12}:\CI_{A,B}^C \cong \CI_{B,A}^C$. 
We rewrite them. 
$$\begin{array}{rl}
\mbox{The left side of (2.3)}=&E(\langle c', \widetilde{\CY}(a,z)e^{L(-1)x}b \rangle)
=E(\langle c', e^{L(-1)x}\widetilde{\CY}(a,z-x)b \rangle) \cr
=&E(\langle e^{L(1)x}c', \widetilde{\CY}(a,z-x)b \rangle) \cr
\mbox{The right side of (2.3)}
=&E(\langle c', \CY(b,x)e^{L(-1)(z)}a \rangle)
=E(\langle c', e^{L(-1)z}\CY(b,x-z)a \rangle) \cr
=&E(\langle e^{L(1)x}c', e^{L(-1)(z-x)}\CY(b,x-z)a \rangle).
\end{array}$$
Since $\langle e^{L(1)x}c', \widetilde{\CY}(a,z-x)b \rangle$ 
and $\langle e^{L(1)x}c', e^{L(-1)z}\CY(b,x-z)a \rangle$ are multivalued 
rational functions on $\{(x,z)\mid x\not=z\}$, 
we may choose $\sigma_{12}$ so that 
$$\sigma_{12}(\CY)(a,z-x)b=e^{L(-1)(z-x)}\CY(b,x-z)a. \eqno{(2.4)}$$  
Similarly, for $\CY\in \CI_{A,B}^{C'}$ and 
canonical intertwining operators $\CY_{C,C'}^{V'}$ and $\CY_{B',B}^{V'}$ 
induced from inner products, there is $\CY^4\in \CI_{A,C}^{B'}$ such that 
$$ E(\langle \1, \CY_{C,C'}^{V'}(c,x)\CY(a,y)b\rangle)=
E(\langle \1, \CY_{B',B}^{V'}(e^{L(-1)(x-y)}\CY^4(a,y-x)c,y)b\rangle).$$
Therefore, we have an isomorphism $\sigma_{23}:\CI_{A,B}^C\cong \CI_{A,C'}^{B'}$. 
We need the notation $\sigma_{23}(\CY)$, but not an explicit formula in this paper.   

In (2.1), we used $\CY^{\boxtimes}$ as the second intertwining operator  
of products. 
Not only the second one, we can also use it for the first 
one at the same time. Actually, for 
$\CY_5(\CY_{A,B}^{\boxtimes})$ with $\CY_5\in \CI_{A\boxtimes B,C}^D$,  
we have $\sigma^{-1}_{123}(\CY_5)\in \CI_{C,D'}^{(A\boxtimes B)'}$ and so 
there is $\delta\in \Hom_V(C\boxtimes D',(A\boxtimes B)')$ 
such that $\sigma^{-1}_{123}(\CY_5)=\delta \CY_{C,D'}^{\boxtimes}$. Therefore we have: 
$$\CY_5(\CY_{A,B}^{\boxtimes})=\sigma_{123}(\delta \CY_{C,D'}^{\boxtimes})(\CY_{A,B}^{\boxtimes})
=\sigma_{123}(\CY_{C,D'}^{\boxtimes})(\delta^{\ast}\CY_{A,B}^{\boxtimes}), $$ 
where $\delta^{\ast}\in \Hom_V(A\boxtimes B, (C\boxtimes D')')$ is a dual of $\delta$ and 
$\sigma_{123}$ denotes $\sigma_{12}\sigma_{23}$.

\subsection{Semi-rigidity and intertwining operators}
We next describe the semi-rigidity in terms of intertwining operators. 
For a $V$-module $U$, let ${\rm rad}^V(U)$ denote the smallest submodule such that 
$U/{\rm rad}^V(U)$ is a direct sum of copies of $V$. 
From the definition of semi-rigidity, if $W$ is not semi-rigid, 
then 
$$\mu(W\boxtimes {\rm rad}^V(\widetilde{W}\boxtimes W))+
\Ker(e_W\boxtimes {\rm id}_W)=
(W\boxtimes \widetilde{W})\boxtimes W, \eqno{(2.5)}$$ 
for any $e_{W}:W\boxtimes \widetilde{W}\to V$, 
where $\mu:(W\boxtimes \widetilde{W})\boxtimes W
\to W\boxtimes (\widetilde{W}\boxtimes W)$ is a canonical isomorphism. 
On the other hand, as we has shown in \S 2.3, 
for any $e_W\CY_{W,\widetilde{W}}^{\boxtimes}\in \CI_{W,\widetilde{W}}^V$, $w, w^1\in W$, 
$\widetilde{w}\in \widetilde{W}$, and $a'\in W'$, 
there is 
$\delta\in \Hom_V(W\boxtimes \widetilde{W}, (W\boxtimes W')')$ such that 
$$
E(\langle a', \sigma_{12}(Y^{W})(w,x)
e_{\widetilde{W}}\CY_{\widetilde{W},W}^{\boxtimes}(\widetilde{w},y)w^1
\rangle)=E(\langle a', \sigma_{123}(\CY_{W,W'}^{\boxtimes})
(\delta \CY_{W,\widetilde{W}}^{\boxtimes}(w,x-y)\widetilde{w},y)
w^1 \rangle). \eqno{(2.6)}$$
Therefore, 
$W$ is not semi-rigid if and only if ${\rm Image}(\delta)$ does not have a factor 
isomorphic to $V$ if and only if 
$\Ker(\delta) +{\rm rad}^V(W\boxtimes \widetilde{W})=W\boxtimes \widetilde{W}$ 
for any $e_{\widetilde{W}}$. 

\subsection{Pseudo-trace}
Although we won't treat pseudo-trace functions \cite{M1} in this paper, 
we will explain them a little. It was introduced to explain 
a symmetric function on $n$-th Zhu algebra $A_n(V)$ in terms 
of $V$-modules.  
Most all of symmetric functions on $A_n(V)$ are 
linear combinations of traces of grade-preserving operators $o(v)$ with $v\in V$ 
on $n$-th lowest homogeneous weight-space $U(n)$ of simple $V$-modules $U$. 
However, in some VOAs, these functions don't cover all 
symmetric functions. The remaining symmetric functions are 
given by the following:
For a $V$-module $U$ with submodules $U\supseteq T\supseteq S$ and 
a surjective $V$-homomorphism $\phi:U \to S$ with $\Ker(\phi)=T$, 
take a transversal $\epsilon:S\to U$, that is, 
$\phi\epsilon=1_S$ and choose a basis $\{s^i \mid i\in I\}$ of $S$. 
We extend $\{s^i, \epsilon(s^i)\}$ to a basis 
$\{s^i, \epsilon(s^i), ... \mid i\in I\}$ of $U$ and using this basis, we can 
express the action of $V$ on $U$ as  
$$ Y^U(v,z)=\begin{pmatrix}A_{11}(v,z)&A_{12}(v,z)&A_{13}(v,z) 
\cr O& A_{22}(v,z)&A_{23}(v,z) \cr O&O& A_{11}(v,z) \end{pmatrix}.\eqno{(2.7)}$$
Then a pseudo-trace function on $(U,\phi)$ is defined by 
$$\Tr^{\phi}_UY^U(z^{L(0)}v,z)q^{\tau(L(0)-c/24)}:=
\sum_{i} \langle  (s^i)', 
Y^U(z^{L(0)}v,z)q^{\tau(L(0)-c/24)}\epsilon(s^i) \rangle, \eqno{(2.8)}$$
where $\{(s^i)'\mid i\in I\}$ is the dual basis of $\{s^i\mid i\in I\}$. 
In other words, it is a trace function of $A_{13}(z^{L(0)}v,z)q^{\tau(L(0)-c/24)}$. 
If it is symmetric with respect to $V$ (we call it $V$-symmetric), that is, 
it is symmetric with the grade-preserving actions of $V$, then 
we call it pseudo-trace (function).  
From (2.7), we have that the value of pseudo-trace 
is zero for an element which acts on $U$ semi-simply. 
For example, $\Tr_U^{\phi}Y^U(\1,z)q^{\tau(L(0)-c/24)}=0$.

\section{Geometrically modified module}
We quote the theory of composition-invertible power series and their actions on modules for 
the Virasoro algebra developed in \cite{H1}. From now on, $q^x$ denotes $e^{2\pi ix}$ for variables $x$ 
to simplify the notation. 
Let $A_j$  $(j=1,2,...)$ be the complex numbers defined by 
$$\frac{1}{2\pi i}(q^y-1)=\left( \exp\left(-\sum_{j=1}^{\infty} A_j y^{j+1}\frac{\pd}{\pd y}\right) \right)y $$
and set 
$$\CU(q^{x})=q^{xL(0)}(2\pi i)^{L(0)}e^{-\sum_{j=1}^{\infty} A_jL(j)}.$$ 
The important one is $\CU(1)$, which satisfies 
$$ \CU(1)\CY(w,x)\CU(1)^{-1}=\CY(\CU(q^{x})w, q^{x}-1)=\CY(q^{xL(0)}\CU(1)w,q^x-1)=\CY[\CU(1)w,x]  \eqno{(3.1)}$$
for an intertwining operator $\CY$, see \cite{Zh} for $\CY[\cdot,x]$.

\subsection{Trace functions}
We first consider $q^{\tau}$-traces of geometrically-modified module 
operators with one more variable $z$:
$$\Psi_U^{\phi}(v;z,\tau):=\Tr^{\phi}_U Y(\CU(q^{z})v, q^{z})q^{\tau(L(0)-c/24)} \eqno{(3.2)}$$
for a $V$-module $U$ and $v\in V$, 
where $\Tr_U^{\phi}$ is a pseudo-trace (including an ordinary trace $\Tr_U$) 
and $c$ is the central charge of $V$. We note that for an ordinary 
trace function, we can consider the trace functions for 
not only $V$ but also a $V$-module $T$ and $\CY\in \CI_{T,U}^U$. 
Namely, we can define a trace function 
$$\Psi_U^{\tr}(\CY;t;z,\tau):
=\Tr_U (\CY(\CU(q^{z})t, q^{z})q^{\tau(L(0)-c/24)}) \qquad t\in T. \eqno{(3.3)}$$

We have to note that $L(0)$ may not be semisimple on a $V$-module $U$. 
We denote the semisimple part of $L(0)$ by $\wt$ and 
$L(0)^{nil}=L(0)-\wt$ is a nilpotent part of $L(0)$. 
Then we will understand $q^{\tau L(0)}$ on $U$ as 
$$q^{\tau L(0)}:=q^{\tau (\wt+L(0)^{nil})}
=q^{\tau \wt}(e^{2\pi i \tau L(0)^{nil}})=
q^{\tau \wt}\sum_{j=0}^{\infty}\frac{(2\pi i\tau L(0)^{nil})^j}{j!}.$$   
In particular, trace function may have a term $q^{\tau r}\tau^j$ for $j\in \N$. 

We note that for simple modules $W$ and $U$, 
$\CY_{W,U}^U\in \CI_{W,U}^U$ has no logarithmic terms and 
the grade-preserving operators $o(w)$ of $w\in W_r$ in 
$\CY_{W,U}^U(w,z)=\sum w_mz^{-m-1}$ is $w_{r-1}$. 
Therefore, by setting $\CU(1)w=\sum w^r$ with homogeneous elements $w^r\in W_r$, we have 
$$\begin{array}{rl}
\Tr^{\phi}_U \CY_{W,U}^U(\CU(q^{z})w, q^{z})q^{\tau(L(0)-c/24)}=&\sum_r \Tr^{\phi}_U 
q^{z(\wt(v^{r}))}w^{r}_{r-1}q^{(-r)}q^{\tau(L(0)-c/24)}\cr
=&\sum_r \Tr^{\phi}_U w^{r}_{r-1} q^{\tau(L(0)-c/24)}.  
\end{array}\eqno{(3.4)}$$ 
Thus, (3.4) is independent of $z$. Moreover, 
it has shown in \cite{H1} that these $q^{\tau}$-traces are absolutely 
convergent when $0<|q^{\tau}|<1$ and can be analytically extended to analytic functions of 
$\tau$ in the upper-half plane.

We next consider $q^{\tau}$-traces of products of two geometrically-modified intertwining operators: 
$$\begin{array}{l}
\Tr^{\phi}_U \CY_1(\CU(q^{y})\CY_{W,\widetilde{W}}^{\boxtimes}(w,x-y)\widetilde{w}, q^{y})q^{\tau(L(0)-c/24)}\cr
\Tr_U^{\phi} \CY_2(\CU(q^{x})w,q^{x})
\CY_{\widetilde{W},U}^{\boxtimes}(\CU(q^{y})\widetilde{w},q^{y})q^{\tau(L(0)-c/24)}
\end{array}\eqno{(3.5)}$$ 
for $w\in W, \widetilde{w}\in \widetilde{W}$, 
$\CY_1\in \CI_{W\boxtimes \widetilde{W},U}^U$, and 
$\CY_2\in \CI_{W,\widetilde{W}\boxtimes U}^U$. 
As we explained, the first function in (3.5) depends on $x-y$, but not on $y$. 
These formal power series (with log-terms) are 
absolutely convergent in $\Omega_1=\{(x,y,\tau)\in \C^2\oplus \CH \mid 0<|q^{x}-q^u|<|q^y|\}$ 
and $\Omega_2=\{(x,y,\tau)\in \C^2\oplus \CH\mid 0<|q^{\tau}|<|q^y|<|q^x|<1\}$, 
respectively, as shown in \cite{H1}, 
where $\CH=\{\tau\in \C\mid {\rm Im}(\tau)>0\}$ is the upper half plane. 
We extend these function analytically to multivalued 
analytic functions on  
$$M_1^2=\{(x,y,\tau)\in \C^2\times \CH \mid x\not=y+p\tau+q \quad\mbox{ for all }
p,q\in \Z\}.$$
We can lift them to 
single valued analytic functions 
$$\begin{array}{rl}
\Psi_U^{\phi}(\CY_1(\CY_{W,\widetilde{W}}^{\boxtimes}):w,\widetilde{w};x,y,\tau):&=E(\Tr^{\phi}_U \CY_1(\CU(q^{y})\CY_{W,\widetilde{W}}^{\boxtimes}(w,x-y)\widetilde{w}, q^{y})q^{\tau(L(0)-c/24)}) \cr
\Psi_U^{\phi}(\CY_2\cdot\CY_{\widetilde{W},U}^{\boxtimes}:w,\widetilde{w};x,y,\tau):&=E(\Tr_U^{\phi} \CY_2(\CU(q^{x})w,q^{x})
\CY_{\widetilde{W},U}^{\boxtimes}(\CU(q^{y})\widetilde{w},q^{y})q^{\tau(L(0)-c/24)})
\end{array}\eqno{(3.6)}$$ 
on the universal covering $\widetilde{M^2_1}$. 
Although Huang has treat only trace functions in \cite{HV}, but it is still possible for 
pseudo-trace functions, (see \cite{M1}).

We need to extend one statement in \cite{HV} 
to logarithmic intertwining operators.

\begin{lmm}
For a (logarithmic) intertwining operator 
$\CY\in \CI_{B,U}^T$, $w\in W$ and $b\in B$, we have 
$$\begin{array}{l}  e^{\tau L(0)}\CY(b,z)u=\CY(e^{\tau L(0)}b,e^{\tau}z)e^{\tau L(0)}u  \qquad 
\cr
q^{\tau L(0)}\CY(\CU(q^{y})b,q^{y})=\CY(\CU(q^{y+\tau})b,q^{y+\tau})q^{\tau L(0)}\qquad \mbox{ and} \cr
\CY^1(\CY^2(\CU(q^y)b,q^y-q^x)\CU(q^x)w,q^x)
=\CY^1(\CU(q^x)\CY^2(b,y-x)w,q^x)
\end{array}$$
\end{lmm}

\pr
Set $\CY(b,z)=\sum_{h=0}^K\sum_{n\in \C} b_{n,h}z^{-n-1}\log^h z$ and 
$y=\log z$ to simplify the notation.
From $\CY(L(-1)b,z)=\frac{d}{dz}\CY(b,z)$, 
we have $(L(-1)b)_{n+1,h}=(-n-1)b_{n,h}+(h+1)b_{n-1,h+1}$ and 
$$L(0)(b_{n,h}u)-b_{n,h}L(0)u\!=\!(L(-1)b)_{n+1,h}+(L(0)b)_{n,h}\!=\!
(-n-1)b_{n,h}\!+\!(h+1)b_{n,h+1}\!+\!(L(0)b)_{n,h}$$
for $u\in U$. Therefore, we obtain:
$$\begin{array}{l}
e^{\tau L(0)}(\sum_{h=0}^K b_{n,h}uy^h)e^{(-n-1)y}
=\!\!\sum_{m=0}^{\infty} \frac{L(0)\tau}{m!}(\sum_h b_{n,h}uy^he^{(-n-1)y}\cr
=\!\!\sum_{m,h,j}\frac{1}{m!}\binom{m}{j}\tau^m 
(L(0)|_b+L(0)|_u\!-\!n\!-\!1)^{m-j}(h\!+\!1)\cdots(h\!+\!j)b_{n,h+j}u(2\pi iy)^he^{(-n-1)y} \cr
=\!\!\sum_{m, k=0}^{\infty}\sum_{j=0}^k
\frac{1}{(m-j)!}\frac{1}{j!}(\tau(L(0)|_b\!+\!L(0)|_u\!-\!n\!-\!1))^{m\!-\!j}
(k\!-\!j\!+\!1)\cdots(k)\tau^j
y^{k-j}b_{n,k}ue^{(-n-1)y}\cr
=\!\!\sum_{k=0}^{\infty}e^{\tau(L(0)|_b+L(0)|_u-(n+1))}b_{n,k}u(y+\tau)^ke^{(-n-1)y}\cr 
=\!\!\sum_k (e^{\tau L(0)}b)_{n,k}(e^{\tau L(0)}u) e^{(-n-1)(y+\tau)}(y+\tau)^k\cr
=\!\!\CY(e^{\tau L(0)}b, e^{\tau+y})e^{\tau L(0)}u=\CY(e^{\tau L(0)}b, e^{\tau}z)e^{\tau L(0)}u, 
\end{array}$$
where $L(0)|_b$ and $L(0)|_u$ denote the action of $L(0)$ on $b$ and $u$, respectively. 
Replacing $\tau$ and $y$ by $2\pi i\tau$ and $2\pi i y$, respectively, 
we have the second equation.
The third comes from $\CU(1)\CY(b,x)=\CY(\CU(q^x)b,q^x-1)\CU(1)$ and the second 
equation. 
\prend

\section{Transformations}
For a $C_2$-cofinite VOA $V$ satisfying Hypothesis I and II, 
$V^{\otimes n}$ is also a $C_2$-cofinite VOA satisfying Hypothesis I and II. 
Moreover, for a $V$-module $W$, $W^{\otimes n}$ is a semi-rigid $V^{\otimes n}$-module 
if and only if $W$ is semi-rigid. We also have that 
$\Psi_V^{\tr}$ appears in $S(\Psi_V^{\tr})$ if and only if 
$\Psi_{V^{\otimes n}}^{\tr}$ appears 
in $S(\Psi_{V^{\otimes n}}^{\tr})$. Therefore, by taking a suitable $V^{\otimes n}$ 
instead of $V$, 
we may assume that $W$ and $\widetilde{W}$ have integer weights to simplify  
the arguments. 

\subsection{Three transformations}
A (pseudo-)trace function of 
$\CY^1(\CY^2)\in \CI_{E,U}^U(\CI_{\widetilde{W},W}^E)$ on $U$ is 
$$\Psi_U^{\phi}(\CY^1(\CY^2):\widetilde{w},w;x,y,\tau)
=E(\Tr^{\phi}_U \CY^1(\CU(q^{y})\CY^2(\widetilde{w},x-y)w,q^{y})q^{\tau(L(0)-c/24)}),
\eqno{(4.1)} $$
for $w\in W,\widetilde{w}\in \widetilde{W}$. A modular transformation $S:\tau\to -1/\tau$ on $\Psi_U^{\phi}$ is 
defined by 
$$\begin{array}{l}
S\left(\Psi_U^{\phi}\right)\left(\CY^1(\CY^2):\widetilde{w},w;x,y,\tau\right)\cr
\mbox{}\qquad=\Psi_U^{\phi}\left(\CY^1(\CY^2): \left( \frac{-1}{\tau}\right)^{L(0)}\widetilde{w}, \left(\frac{-1}{\tau}\right)^{L(0)}w;
\frac{-1}{\tau}x,\frac{-1}{\tau}y; \frac{-1}{\tau}\right).\end{array}\eqno{(4.2)}$$
When $\CY^1(\CY^2)=Y^U(\CY)$ for some $\CY\in \CI_{\widetilde{W},W}^V$, 
it has a modular invariance property. In other words, there 
are $\lambda_{(T,\psi)}\in \C$ such that 
$$
S\left(\Psi_U^{\phi}\right)\left(Y^U(\CY):\widetilde{w},w;x,y,\tau \right)
=\sum \lambda_{(T,\psi)} \Psi_T^{\psi}\left( Y^T(\CY):\widetilde{w},w;x,y,\tau
\right).\eqno{(4.3)}$$
We note that $\lambda_{(T,\psi)}$ does not depend on 
$\CY\in \CI_{W,\tilde{W}}^V$, but on $V$. 

We define actions $S$, $\alpha_t$, $\beta_t$ on $R_2^1$ by 
$$ \begin{array}{ccc}
(x,y,\tau) &\xrightarrow{S} & (-x/\tau,-y/\tau, -1/\tau) \cr
\downarrow \beta_t & & \downarrow \alpha_t \cr
(x,y+t,\tau) &\xrightarrow{S} & (-x/\tau, -y/\tau+1,-1/\tau). 
\end{array}\eqno{(4.4)}$$
Along a line $\CL=\{(x,y+t,\tau)\mid t\in [0,1]\}$ from 
$(x,y,\tau)$ to $(x,y+1,\tau)$, we define 
$$
\alpha_t(\Psi_U^{\phi})(\CY:\widetilde{w},w;x,y,\tau):
=\Psi_U^{\phi}(\CY:\widetilde{w},w;x,y+t,\tau).
\eqno{(4.5)}$$
Since $(x,y,\tau)\to (x,y+t,\tau)$ preserves 
$\Omega_2=\{(x,y,\tau)\in \C^2\oplus H\mid |q^{\tau}|<|q^y|<|q^x|<1\}$, we have 
$$\begin{array}{rl}
\alpha_t(\Psi_U^{\phi})(\CY^1(\CY^2):\widetilde{w},w;x,y,\tau)
=&\alpha_t(\Tr_U^{\phi} \CY_3(\CU(q^{x})\widetilde{w},q^{x})
\CY_{\widetilde{W},U}^{\boxtimes}(\CU(q^{y})w,q^{y})q^{\tau(L(0)-c/24)})\cr
=&\Tr_U^{\phi} \CY^3(\CU(q^{x})\widetilde{w},q^{x})
\CY_{\widetilde{W},U}^{\boxtimes}(\CU(q^{y+t})w,q^{y+t})q^{\tau(L(0)-c/24)}\cr
=&\Tr^{\phi}_U \CY^4(\CU(q^{y})\CY^5(\widetilde{w},x-y)w, q^{y})q^{\tau(L(0)-c/24)}\cr
=&\Psi_U^{\phi}(\CY^4(\CY^5):\widetilde{w},w;x,y,\tau)
\end{array}\eqno{(4.6)}$$
for some $\CY^3$ and $\CY^4(\CY^5)$, because $\CY_{W,U}^{\boxtimes}(\CU(q^{y+t})w,q^{y+t})$ 
is a linear combination of geometrically modified intertwining 
operators in $\CI_{W,U}^{\boxtimes}$. 

An important case is where $U=V$ and $\CY^1(\CY^2)=Y(\CY)$ 
with $\CY\in \CI_{\widetilde{W},W}^V$. 
Then since $W\boxtimes V=W$ is irreducible, 
$$\alpha_1(\Psi_V^{\tr})(Y(\CY))=e^{2\pi i\wt(W)}\Psi_V^{\tr}(Y(\CY)).$$ 
We set $\kappa=e^{2\pi i\wt(W)}$.  
We then define $\beta_t$ according to a line $S^{-1}(\CL)$ by 
$$\beta_t(\Psi_U^{\phi})(\CY^1(\CY^2):\widetilde{w},w;x,y,\tau)
=\Psi_U^{\phi}(\CY^1(\CY^2):\widetilde{w},w;x,y+t\tau,\tau) \mbox{ for any }\Psi_U^{\phi}.  \eqno{(4.7)}$$
Since $\alpha S=S\beta$ on $R_2^1$ and 
$$S(\Psi_U^{\phi})(\CY^1(\CY^2):\widetilde{w},w)=(-1/\tau)^{(\wt(w)+\wt(\widetilde{w}))}
\Psi_U^{\phi}(\CY^1(\CY^2):\widetilde{w},w)S,$$ 
we have the following relation.

\begin{prn}  
$$\beta_t (S(\Psi_V)))=S(\alpha_t(\Psi_V)).  \eqno{(4.8)}$$ 
\end{prn}

For $\CY\in \CI_{W,\tilde{W}}^V$, we have  
$$S(\Psi_V^{\tr})(Y(\CY))=\sum_{U} \lambda_{(U,\tr)}\Psi_U^{\tr}(Y(\CY))$$
by Hypothesis II and we will consider the following diagram: 
$$\begin{array}{cccc}
    \Psi_V^{\tr}(Y(\CY)) & \xrightarrow{\alpha} &  &\kappa \Psi_V^{\tr}(Y(\CY)) \cr
      \downarrow \quad S        &   & & \downarrow \quad S  \cr
    \sum \lambda_{(U,\tr)}\Psi_U^{\tr}(Y^U(\CY)) &\xrightarrow{\beta}  & 
\sum \lambda_{(U,\tr)}\beta(\Psi_U^{\tr}(Y^U(\CY))&=\kappa\sum \lambda_{(U,\tr)}\Psi_U^{\tr}(Y^U(\CY)) \cr
    \end{array}$$

\subsection{The image of $\beta$}
In this section, we will calculate 
$\beta_1(\Psi_{U}^{\tr})(Y^U(\CY_{\widetilde{W},W}^V))$ as a formal power series.  
In other words, we expand them in the area  
$0<|q^y|<|q^x|$ and $0<|q^{\tau}|<1$ as formal (rational) power series 
of $(x-y)$ and $q^{\tau}$ and $\tau$. We note $|q^{y+t\tau}|\leq |q^y|<|q^x|$. 

Set $A=(W\boxtimes U)$ and $\CY_{\widetilde{W},W}^V
=e_{\widetilde{W}}\CY_{\widetilde{W},W}^{\boxtimes}$, then we have:
$$\begin{array}{l}
\beta_1(\Psi^{\tr}_U)(Y^U(\CY_{\widetilde{W},W}^V):\widetilde{w},w;x,y,\tau) \cr 
\mbox{}=E(\Tr_U Y^{U}(\CU(q^{y+\tau})
\CY_{\widetilde{W},W}^V(\widetilde{w},x-(y+\tau))w, q^{y+\tau})q^{\tau(L(0)-\frac{c}{24})})
\quad \cr
\mbox{}=E(\Tr_U Y^{U}(\CY_{\widetilde{W},W}^V(\CU(q^x)\widetilde{w},q^{x}-q^{y+\tau})
\CU(q^{y+\tau})w, 
q^{y+\tau})q^{\tau(L(0)-\frac{c}{24})})\quad \mbox{ by Lemma 1}\cr
\mbox{}=E(\Tr_{U} \sigma_{23}(\CY_{\widetilde{W},U'}^{\boxtimes})(\CU(q^x)\widetilde{w},q^x)\xi_U
\CY_{W,U}^{\boxtimes}
(\CU(q^{y+\tau})w,q^{y+\tau})q^{\tau(L(0)-\frac{c}{24})})\cr
\mbox{}\qquad \qquad \mbox{for some }\xi_U\in \Hom_V(W\boxtimes U,(\widetilde{W}\boxtimes U')') \cr
\mbox{}=E(\Tr_U \sigma_{23}(\CY_{\widetilde{W},U'}^{\boxtimes})(\CU(q^x)\widetilde{w},q^x)
q^{\tau(L(0)-\frac{c}{24})}
\xi_U\CY_{W,U}^{\boxtimes}(\CU(q^{y})w,q^y))
\qquad \mbox{ by Lemma 1}\cr
\mbox{}=E(\Tr_U \sigma_{23}(\CY_{\widetilde{W},U'}^{\boxtimes})(\CU(q^x)\widetilde{w},q^x)
\xi_U q^{\tau(L(0)-\frac{c}{24})}
\CY_{W,U}^{\boxtimes}(\CU(q^{y})w,q^y))\cr
\mbox{}=E(\Tr_A \CY_{W,U}^{\boxtimes}(\CU(q^y)w,q^y)
\sigma_{23}(\CY_{\widetilde{W},U'}^{\boxtimes})(\CU(q^x)\widetilde{w},q^x)\xi_U q^{\tau(L(0)-c/24)})\cr
\mbox{}\qquad \qquad \mbox{because the trace is symmetric}\cr
\mbox{}=E(\Tr_A \sigma_{123}(\CY_{A,A'}^{\boxtimes}(\delta_U \CY_{W,\widetilde{W}}^{\boxtimes})(\CU(q^y)w,
q^y-q^x)\CU(q^x)\widetilde{w},q^x)q^{\tau(L(0)-c/24)}) \cr
\mbox{}\qquad \qquad \mbox{for some $\delta_U\in \Hom_V(W\boxtimes \widetilde{W},(A\boxtimes A')')$.} 
\end{array}$$
Set $L[-1]=L(-1)+L(0)$ (see \cite{Zh}). 
Then we get $U(1)e^{L(-1)z}=e^{(2\pi i)L[-1]z}U(1)$ from (3.1) 
and so the above equals to the following:
$$\begin{array}{l}
E(\Tr_A \sigma_{123}(\CY_{A,A'}^{\boxtimes})(\delta_U 
\CU(q^x)\CY_{W,\widetilde{W}}^{\boxtimes}(w,y-x)\widetilde{w},q^x)q^{\tau(L(0)-c/24)}) 
\qquad\mbox{ by Lemma 1}\cr
\mbox{}=E(\Tr_A \sigma_{123}(\CY_{A,A'}^{\boxtimes})(\delta_U 
q^{L(0)x}\CU(1)e^{L(-1)(y-x)}\sigma_{12}(\CY_{W,\widetilde{W}}^{\boxtimes})
(\widetilde{w}, x-y)w,q^x)q^{\tau(L(0)-c/24)}) \cr
\mbox{}\qquad \qquad \mbox{by skew symmetry intertwining operator, see (2.4)} \cr
\mbox{}=E(\Tr_A \sigma_{123}(\CY_{A,A'}^{\boxtimes})(\delta_U 
q^{L(0)x}e^{(2\pi i)L[-1](y-x)}\CU(1)\sigma_{12}(\CY_{W,\widetilde{W}}^{\boxtimes})
(\widetilde{w}, x-y)w,q^x)q^{\tau(L(0)-c/24)}).
\end{array}$$
As we explained, the pair of terms $q^{L(0)x}$ and $q^x$ in the 
above expression is just formal and has 
no influence. The next term is $e^{(2\pi i)L[-1](y-x)}$. However, since 
the grade preserving operators of $L[-1]u$ are zero for any 
$u\in \widetilde{W}\boxtimes W$, we finally have 
$$\begin{array}{l}
\beta_1(\Psi^{\tr}_U)(Y^U(\CY_{\widetilde{W},W}^V):\widetilde{w},w;x,y,\tau) \cr 
\mbox{}\quad=E(\Tr_A \sigma_{123}(\CY_{A,A'}^{\boxtimes})
\CU(q^x)\delta_U\sigma_{12}(\CY_{W,\widetilde{W}}^{\boxtimes}))
(\widetilde{w}, x-y)w,q^x)q^{\tau(L(0)-c/24)}). 
\end{array}\eqno{(4.9)}$$
In particular, we have the following lemma. \\

\noindent
{\bf Lemma 3}\quad {\it
$\beta_1(\Psi_{U}^{\tr})(Y^U(\CY_{\widetilde{W},W}^V))$ is again 
an ordinary trace function.}\\

We express the definitions of $\xi_U$ and $\delta_U$ in a short way 
$$
Y^U(\CY_{\widetilde{W},W}^V)
=\sigma_{23}(\CY_{\widetilde{W},U'}^{\boxtimes})\xi\CY_{W,U}^{\boxtimes}
\quad\mbox{ and }\quad 
\CY_{W,U}^{\boxtimes}\sigma_{23}(\CY_{\widetilde{W},U'}^{\boxtimes})\xi_U=
\sigma_{123}(\CY_{A,A'}^{\boxtimes})(\delta_U \CY_{W,\widetilde{W}}^{\boxtimes}).
\eqno{(4.10)}$$ 

For $a'\in A'$, $\widetilde{w}\in \widetilde{W}$, $w,w^1\in 
W$ and $u\in U$, let us consider 
$$\begin{array}{l}
\langle a', \CY_{W,U}^{\boxtimes}(w^1,x)
Y^U(e_{\widetilde{W}}\CY_{\widetilde{W},W}^{\boxtimes}(\widetilde{w},y-z)w,z)u\rangle 
\end{array} \eqno{(4.11)}$$
into two ways. 
Set $B={\rm Image}(\delta_U)$, then 
there is $\CY_{B,W}^{(U\boxtimes A')'}$ such that  
$$\begin{array}{rl}
\mbox{(4.11)}&=
\langle a', \CY_{W,U}^{\boxtimes}(w^1,x)
\sigma_{23}(\CY_{\widetilde{W},U'}^{\boxtimes})(\widetilde{w},y)
\xi\CY_{W,U}^{\boxtimes}(w,z)u\rangle \cr
&=
\langle a', \sigma_{123}(\CY_{A,A'}^{\boxtimes})
(\delta_U\CY_{W,\widetilde{W}}^{\boxtimes}(w^1,x-y)\widetilde{w},y)
\CY_{W,U}^{\boxtimes}(w,z)u\rangle \cr
&=
\langle a', \sigma_{123}(\CY_{U,A'}^{\boxtimes})
\CY_{B,W}^{(U\boxtimes A')'}(\delta_U\CY_{W,\widetilde{W}}^{\boxtimes}(w^1,x-y)\widetilde{w},y-z)
w,z)u\rangle. 
\end{array}$$
On the other hand, there is $\CY_{W,V}^{W}
\in \CI_{W,V}^{W}$ and $\epsilon\in \Hom_V(W,(U\boxtimes A')')$ such that 
$$\begin{array}{rl}
\mbox{(4.11)}=& 
\langle a', \CY_{W,U}^{\boxtimes}(Y_{W,V}^{W}
(w^1,x-z)
e_{\widetilde{W}}\CY_{\widetilde{W},W}^{\boxtimes}(\widetilde{w},y-z)w,z)u\rangle \cr
=&\langle a', \sigma_{123}(\CY_{U,A'}^{\boxtimes})(\epsilon Y_{W,V}^{W}
(w^1,x-z)
e_{\widetilde{W}}\CY_{\widetilde{W},W}^{\boxtimes}(\widetilde{w},y-z)w,z)u\rangle  
\end{array}$$
for any $a'\in A'$ and $u\in U$. We note $\CY_{W,V}^W\in \C\sigma_{12}(Y^W)$. 
Therefore, we have 
$$
\epsilon Y_{W,V}^{W}
(w^1,x-z)
e_{\widetilde{W}}\CY_{\widetilde{W},W}^{\boxtimes}(\widetilde{w},y-z)w
=\CY_{B,W}^{(U\boxtimes A')'}
(\delta_U\CY_{W,\widetilde{W}}^{\boxtimes}(w^1,x-z)\widetilde{w},y-z)w.$$
Since the image of $\epsilon$ is $W$, we obtain 
$$
\epsilon Y_{W,V}^{W}
(w^1,x-z)
\CY_{\widetilde{W},W}^{\boxtimes}(\widetilde{w},y-z)w
=\CY_{B,W}^W
(\delta_U\CY_{W,\widetilde{W}}^{\boxtimes}(w^1,x-z)\widetilde{w},y-z)
w$$
for some $\CY_{B,W}^W$. 
Thus, $\delta_U$ in (4.8) essentially coincides with $\delta$ in (2.6), which does not 
depend on the choice of $U$.

\section{Proof of the Main Theorem}
We now start the proof of the Main Theorem.  
Let $W$ be an irreducible module. 
As we showed in the previous section, 
$$\begin{array}{rl}
 \beta_1(\sum \lambda_{(U,\tr)}\Psi_U^{\tr})(Y(e_{\widetilde{W}}\CY_{\widetilde{W},W}^{\boxtimes}))
 =&\sum \lambda_{(U,\tr)}\beta_1(\Psi_U^{\tr})(Y(e_{\widetilde{W}}\CY_{\widetilde{W},W}^{\boxtimes})) \cr
 =&\sum \lambda_{(U,\tr)}\Psi_{W\boxtimes U}
(\CY_{B,U}^U(\delta \CY_{W,\widetilde{W}}^{\boxtimes})). 
\end{array}$$
On the other hand, since $\beta_1(S(\Psi_V))=S(\alpha_1(\Psi_V))$, we obtain 
$$\beta_1(\sum \lambda_{(U,\tr)}\Psi_U^{\tr}(Y(e_{\widetilde{W}}\CY_{\widetilde{W},W}^{\boxtimes}))
=\kappa(\sum \lambda_{(U,\tr)}\Psi_U^{\tr}(Y(e_{\widetilde{W}}\CY_{\widetilde{W},W}^{\boxtimes})).$$
Therefore, we have 
$$\sum \lambda_{(U,\tr)}\Psi_{W\boxtimes U}
(\CY_{B,U}^U(\delta \CY_{W,\widetilde{W}}^{\boxtimes}))
=\kappa(\sum \lambda_{(U,\tr)}\Psi_U^{\tr}(Y(e_{\widetilde{W}}\CY_{\widetilde{W},W}^{\boxtimes})).$$
Suppose that $W$ is not semi-rigid. 
As we mentioned, we may assume that a conformal weight $\wt(W)$ of $W$ is an integer. 
Then $\Ker(\delta)+\Ker(e_{\widetilde{W}})=\widetilde{W}\boxtimes W$. 
Set $Q=\Ker(\delta)\cap \Ker(e_{\widetilde{W}})$ and $W\boxtimes \widetilde{W}/Q=Q^1\oplus Q^2$ with 
$Q^1=\Ker(e_{\widetilde{W}})/Q$ and $Q^2=\Ker(\delta)/Q\cong V$. 
Then $\Psi_{W\boxtimes U}
(\CY_{B,U}^U(\delta \CY_{W,\widetilde{W}}^{\boxtimes}))$ are all given by 
traces on $Q^1$ and $\Psi_U^{\tr}(Y(e_{\widetilde{W}}\CY_{\widetilde{W},W}^{\boxtimes})$ are 
all given by traces on $Q^2$. We hence have 
$$
 \sum \lambda_{(U,\tr)}\Psi_{W\boxtimes U}
(\CY_{B,U}^U(\delta \CY_{W,\widetilde{W}}^{\boxtimes}))=0, $$
which contradicts to 
$\sum \lambda_{(U,\tr)}\Psi_U^{\tr}(Y(e_{\widetilde{W}}\CY_{\widetilde{W},W}^{\boxtimes})
\not=0$. 
Therefore, $W$ is semi-rigid. Since $W$ is arbitrary, $V$ is semi-rigid. \\

We next show $\lambda_{(V',\tr)}\not=0$. Choose a simple module $U$ so that 
$\lambda_{(U,\tr)}\not=0$. Set $W=U'$ and consider the trace function of 
the $e_{\widetilde{W}}\CY_{\widetilde{W},W}^{\boxtimes}$ in 
$\beta_1(\Psi_U^{\tr})(e_{\widetilde{W}}\CY_{\widetilde{W},W}^{\boxtimes})$. 
It has a nonzero scalar multiple of 
$$\Psi_{W\boxtimes U}^{\tr}(e_{\widetilde{W}}\CY_{\widetilde{W},W}^{\boxtimes})$$
and so it has a term 
$\Psi_{V'}^{\tr}(e_{\widetilde{W}}\CY_{\widetilde{W},W}^{\boxtimes})$ with 
a nonzero coefficient.  
On the other hand, for any $V$-modules $T\not=U$, 
$\beta_1(\Psi_T^{\tr}(e_{\widetilde{W}}\CY_{\widetilde{W},W}^{\boxtimes})$ has no entries 
of $\Psi_{V'}^{\tr}(e_{\widetilde{W}}\CY_{\widetilde{W},W}^{\boxtimes})$. 
Therefore, $\Psi_{V'}^{\tr}(e_{\widetilde{W}}\CY_{\widetilde{W},W}^{\boxtimes})$ has 
nonzero coefficient in 
$\beta_1(\sum \lambda_{(U,\tr)}
\Psi_U^{\tr}(Y(e_{\widetilde{W}}\CY_{\widetilde{W},W}^{\boxtimes}))$.

The remaining thing is to prove $\lambda_{(U,\tr)}\not=0$ for 
every simple module $U$. 
Set $W=U'$.
As we showed, $\lambda_{(V',\Tr)}\not=0$ and so 
there is a simple $V$-module $T$ with $\lambda_{(T,\tr)}\not=0$ such that 
$\beta_1(\Psi_T^{\phi})(Y^T(\CY_{\widetilde{W},W}^V)$ to have nonzero coefficient 
at $\Psi_{V'}^{\Tr}(Y^U(\CY_{\widetilde{W},W}^V)$. Then since 
$\Hom_V(T\boxtimes W, V')\not=0$, $T=(W)'=U$ and so $\lambda_{(U,\tr)}\not=0$ as 
we desired. \\
This completes the proof of the Main theorem. \\

\section{Orbifold model}
At last, we will show an example satisfying Hypothesis II. 
Let $T$ be a rational vertex operator algebra of CFT type and 
$\tau\in {\rm Aut}(T)$ of order $p$.  
Let $\xi\in \C$ be a primitive $p$-th root of unity and 
decompose $T$ into $T=\oplus_{i=0}^{p-1}T^{(i)}$ 
with $T^{(i)}=\{v\in T\mid \sigma(v)=\xi^iv\}$.  
We assume that the fixed point subVOA $V:=T^{\tau}$ is $C_2$-cofinite 
and satisfies Hypothesis I. \\

\noindent
{\bf Theorem 4} \quad {\it  
Let $T$ be a rational vertex operator algebra of CFT type and 
$\tau$ is a finite automorphism of $T$.  We assume that 
the fixed point subVOA $T^{\tau}$ is $C_2$-cofinite and satisfies Hypothesis I. 
Then $S(\Psi_{T^{\tau}}^{\tr})$ is a linear combination of trace functions. }\\ 

Before we start the proof of Theorem 4, we first show the following:\\

\noindent
{\bf Proposition 5} \quad {\it  
Under the assumption in Theorem 4, 
$T^{\tau}$ is projective as a $T^{\tau}$-module. }\\

\pr
Suppose false and let $0\to B\xrightarrow{\epsilon} P\xrightarrow{\phi} T^{\tau}\to 0$ 
be a non-split extension of $T^{\tau}$.  Set $V=T^{\tau}$. 
Viewing $T$ as a $V$-module, we define a fusion product $W=T\boxtimes_V P$ and 
set $W^{(i)}=T^{(i)}\boxtimes_V P$. We note 
$W=W^{(0)}\oplus \cdots \oplus W^{(n-1)}$ and $W^{(0)}=P$. Similarly, 
we set $R=({\rm id}_T\boxtimes \epsilon)(T\boxtimes_V B)\subseteq T\boxtimes_VP$ and 
$R^{(i)}=({\rm id}_{T^{(i)}}\boxtimes \epsilon)(T^{(i)}\boxtimes_V B)\subseteq T^{(i)}\boxtimes_VP$.  
We note that $({\rm id}_{T^{(i)}}\boxtimes \epsilon$ may not be injective, but 
$R^i$ is not zero since there is a canonical epimorphism 
$\widetilde{T^{(i)}}\boxtimes (T^{(i)}\boxtimes P)
\cong (\widetilde{T^{(i)}}\boxtimes T^{(i)})\boxtimes P \to V\boxtimes P\cong P$.

As we explained, there is $\CY\in \CI_{T,W}^W$ such that 
$$ E(\langle w', \CY(t,z_1)\CY_{T,P}^{\boxtimes}(t^1,z_2)p\rangle) 
=E(\langle w', \CY_{T,P}^{\boxtimes}(Y(t,z_1-z_2)t^1,z_2)p\rangle)$$
for $t,t^1\in T$, $w'\in W'$ and $p\in P$.  
From the definition of $\CY$ and the Commutativity of vertex operators of $T$, we have 
$$ \begin{array}{rl}
E(\langle w',\CY(t^1,z_1)\CY(t^2,z_2)\CY_{T,P}^{\boxtimes}(t^1,z)p\rangle) 
=&\!E(\langle w',\CY(t^1,z_1)\CY_{T,P}^{\boxtimes}(Y(t^2,z_2-z)t^3,z)p\rangle)\cr 
=&\!E(\langle w',\CY_{T,P}^{\boxtimes}(Y(t^1,z_1-z)Y(t^2,z_2-z)t^3,z)p\rangle)\cr
=&\!E(\langle w',\CY_{T,P}^{\boxtimes}(Y(t^2,z_2-z)Y(t^1,z_1-z)t^3,z)p\rangle)\cr
=&\!E(\langle w',\CY(v^2,z_2)\CY(t^1,z_1)\CY_{T,P}^{\boxtimes}(t^3,z)p\rangle)
\end{array}$$
for $t^1,t^2,t^3\in T$, 
which implies the Commutativity of $\{\CY(t,z) \mid t\in T\}$. We also have 
$$ \begin{array}{rl}
E(\langle w',\CY(t^1,z_1)\CY(t^2,z_2)\CY_{T,P}^{\boxtimes}(t^3,z)p\rangle) 
=&\!E(\langle w',\CY_{T,P}^{\boxtimes}(Y(t^1,z_1-z)Y(t^2,z_2-z)t^3,z)p\rangle)\cr
=&\!E(\langle w',\CY_{T,P}^{\boxtimes}(Y(Y(t^1,z_1-z_2)t^2,z_2-z)t^3,z)p\rangle)\cr
=&\!E(\langle w',\CY(Y(t^1,z_1-z_2)t^2,z_2)\CY_{T,P}^{\boxtimes}(t^3,z)p\rangle).
\end{array}$$
Furthermore, taking $t^1=\1$, we obtain $\CY(t,z)p=\CY_{T,P}^{\boxtimes}(t,z)p$ for 
$t\in V, p\in P$ since
$$\begin{array}{rl}
E(\langle w',\CY(t,z_1)p\rangle)
=&E(\langle w',\CY(t,z_1)\CY_{T,P}^{\boxtimes}(\1,z_2)p\rangle)
=E(\langle w',\CY_{T,P}^{\boxtimes}(Y(t,z_1-z_2)\1,z_2)p\rangle)\cr
=&E(\langle w',\CY_{T,P}^{\boxtimes}(e^{(z_1-z_2)L(-1)}t,z_2)p\rangle) 
=E(\langle w',\CY_{T,P}^{\boxtimes}(t,z_2+z_1-z_2)p\rangle)\cr
=&E(\langle w',\CY_{T,P}^{\boxtimes}(t,z_1)p\rangle). 
\end{array}$$
Therefore, $T\boxtimes_VP$ is a $T$-module and 
$({\rm id}_V\boxtimes \epsilon)(T\boxtimes_VB)$ is a direct summand of 
$T\boxtimes_VP$ since $T$ is rational. 
Then $B=({\rm id}_V\boxtimes \epsilon)(T\boxtimes_VB)\cap T^{(0)}_VP$ is 
also a direct summand of $P$ as a $V$-module, which contradicts the choice of $P$. 
\prend

We will assert one more general result. \\

\noindent
{\bf Proposition 6}\quad {\it
$T^{(i)}$ is a simple current as a $V$-module, that is, 
$T^{(i)}\boxtimes_V D$ is simple for any simple $V$-module $D$. }\\

\pr 
Set $Q=T\boxtimes_V D$ and $Q^{(i)}=T^{(i)}\boxtimes_V D$. 
For $Q$, we will use the same arguments as above. 
Suppose that $Q^{(i)}$ contains a proper submodule $S$. 
Then $S^{\perp}\cap (Q^{(i)})'\not=0$ and so we have 
$$\begin{array}{rl}
E(\langle d',\CY(t^{(i)},z_1)\CY(t^{(n-i)},z)s\rangle)=&E(\langle d',\CY(Y(t^{(i)},z_1-z)t^{(n-i)},z)s\rangle) \cr
=&E(\langle d',Y(Y(t^{(i)},z_1-z)t^{(n-i)},z)s\rangle)=0
\end{array}$$ 
for $t^{(i)}\in T^{(i)}$, $d'\in S^{\perp}$ and $s\in S$ since $Y(t^{(i)},z)t^{(n-i)}\in V\{z\}[\log z]$.  
On the other hand, since 
$E(\langle q',\CY(t^{(i)},z_1)\CY(t^{(n-i)},z)s\rangle)=
E(\langle q',\CY(Y(t^{(i)},z_1-z)t^{(n-i)},z)s\rangle)\not=0$ 
for some $q'\in (Q^{(i)})', t^{(i)}\in T^{(i)}$ and $t^{(n-i)}\in T^{(n-i)}$, 
the coefficients in $\{\CY(t^{(n-i)},z)s\mid s\in S, t^{(n-i)}\in T^{(n-i)}\}$ 
spans $D$ and so those in 
$\{ \CY(t^{(i)},z_1)\CY(t^{(n-i)},z)s\mid t^{(i)}\in T^{(i)}, t^{(n-i)}\in T^{(n-i)} \}$ 
spans $Q^{(i)}$. Therefore, we have a contradiction. \prend

We now start the proof of Theorem 4. 
We pick up one twisted simple $T$-module $M=\oplus_{n=0}^{\infty}M_{\lambda+n/p}$. 
Then for each $i$, 
$W^{(i)}=\oplus_{n=0}^{\infty}M_{\lambda+n+i/p}$ is a simple $V$-module and 
we may assume that $T^{(j)}\boxtimes W^{(i)}=W^{(i+j)}$ since $T^{(j)}$ is simple 
current. 
Using $W=W^{(0)}$ and $\widetilde{W}$, 
we will consider geometrically modified trace functions. 
Set $\CY=\CY_{\widetilde{W},W}^V$. 

Let us consider the images of $\Psi_T(Y(\CY))$ by $\alpha_1$ and $S$. 
Since $W\boxtimes T^{(i)}=W^{(i)}$, 
$\alpha_1(\Psi_{T^{(i)}}^{\tr})(Y(\CY))=e^{2\pi i\wt(W^{(i)})}
\Psi_{T^{(i)}}^{\tr}(Y(\CY))$ by (4.6). Therefore, we have:
$$\alpha_1(\Psi_T^{\tr}(Y(\CY)))=\alpha_1(\sum_{i=0}^{p-1} 
\Psi_{T^{(i)}}^{\tr}(Y(\CY)))
=e^{2\pi i\wt(W^{(0)})}(\sum_{i=0}^{p-1} \xi^i\Psi_{T^{(i)}}^{\tr}(Y(\CY))),$$ 
which coincides with a scalar multiple of a $\tau$-twisted trace function  
$$\Psi_T^{\tr}(\tau\cdot Y(\CY):w,\widetilde{w};x,y,\tau):=
E(\Tr_T \tau Y^T(\CU(q^{y})
\CY_{W,\widetilde{W}}^V(w,x-y)\widetilde{w}, q^{y})q^{\tau(L(0)-c/24)})$$
on $T$ with an action of $\tau$. 
On the other hand, 
since $T$ is rational and $C_2$-cofinite, 
$S(\Psi_T^{\tr})$ is a linear combination of trace functions 
$\Psi_U^{\tr}$ on $T$-modules $U$ which is also a $V$-module.   
Therefore, $\beta_1(S(\Psi_T^{\tr}(Y(\CY)))
=S(\alpha_1(\Psi_T^{\tr}(Y(\CY))))=e^{2\pi i\wt(W)}
S(\Psi_T^{\tr}(\tau\cdot Y(\CY)))$ is also a linear combination 
of trace functions. 
Since $\Psi_V^{\tr}=\frac{1}{p}(\sum_{i=0}^{p-1} \Psi_{T}^{\tr}(\tau^i Y(\CY))))$, 
we have the desired conclusion. \\
This completes the proof of Theorem 4. \\

Let us go back to the assumptions in Theorem 4. 
Since $S(\Psi_V^{\tr})$ is a linear combination of trace functions, 
$V$ satisfies the conditions of the main theorem and so $V$ is semi-rigid. 
We have also proved that $V$ is projective as a $V$-module. 
Therefore, we have the following by Corollary 15 in \cite{M1}. \\

\noindent
{\bf Corollary 7}\quad {\it 
Under the assumptions in Theorem 4, $T^{\tau}$ is rational. }\\

\end{document}